\documentclass[12pt]{article}
\usepackage{amsmath,amssymb,amsfonts}     % AMS LateX is standard nowadays.
\usepackage{graphics}                     % For picture environment.
\usepackage{graphicx}                     % For picture environment.
\date{}
\newtheorem{lemma}{Lemma}[section]
\newtheorem{theorem}{Theorem}[section]

\newcommand{\adfQED}{\hfill $\square$}

\newcommand{\adfmod}[1]{~(\mathrm{mod}~#1)}

\newcommand{\adfhide}[1]{}

\newcommand{\adfPgap}{\vskip 2mm}

\newenvironment{adfenumerate}{
\begin{enumerate}
\setlength{\itemsep}{0.5mm}
\setlength{\parskip}{0mm}
\setlength{\parsep}{0mm}
}{
\end{enumerate}
}

\begin{document}
\title{\bf Group divisible designs with block size four and type $g^u b^1 (gu/2)^1$\\}
\author{A.~D.~Forbes\\
        School of Mathematics and Statistics,\\
        The Open University,\\
        Walton Hall, Milton Keynes MK7 6AA, UK
%       \\ Version 3.14, 5 Jun 2019  % HIDE FOR FINAL SUBMISSION.
        }
\maketitle
\begin{abstract}
We discuss group divisible designs with block size four and type $g^u b^1 (gu/2)^1$, where
$u = 5$, 6 and 7. For integers $a$ and $b$, we prove the following.
  (i) A 4-GDD of type $(4a)^5 b^1 (10a)^1$ exists if and only if $a \ge 1$, $b \equiv a$ (mod 3) and $4a \le b \le 10a$.
  (ii) A 4-GDD of type $(6a+3)^6 b^1 (18a+9)^1$ exists if and only if $a \ge 0$, $b \equiv 3$ (mod 6) and $6a+3 \le b \le 18a + 9$.
  (iii) A 4-GDD of type $(6a)^6 b^1 (18a)^1$ exists if and only if $a \ge 1$, $b \equiv 0$ (mod 3) and $6a \le b \le 18a$.
  (iv) A 4-GDD of type $(12a)^7 b^1 (42a)^1$ exists if and only if $a \ge 1$, $b \equiv 0$ (mod 3) and $12a \le b \le 42a$,
  except possibly
  for $12a \in \{120, 180, 240, 360, 420, 720, 840\}$, $24a < b < 42a$,
  for $12a \in \{144, 1008\}$, $30a < b < 42a$, and
  for $12a \in \{168, 252, 336, 504, 1512\}$, $36a < b < 42a$.
\end{abstract}

{\small Keywords: Group divisible design, 4-GDD}

{\small Mathematics Subject Classification: 05B05}

\vskip 10mm

%%%%%%%%%%%%%%%%%%%%%%%%%%%%%%%%%%%%%%%%%%%%%%%%%%%%%%%%%%%%%%%%%%%%%%%%%%%%%%%%%%%%%%%%%%
%%%%%%%%%%%%%%%%%%%%%%%%%%%%%%%%%%%%%%%%%%%%%%%%%%%%%%%%%%%%%%%%%%%%%%%%%%%%%%%%%%%%%%%%%%
%%%%%%%%%%%%%%%%%%%%%%%%%%%%%%%%%%%%%%%%%%%%%%%%%%%%%%%%%%%%%%%%%%%%%%%%%%%%%%%%%%%%%%%%%%

\section{Introduction}\label{sec:Introduction}

A {\em group divisible design}, $K$-GDD, of type $g_1^{u_1} \dots g_r^{u_r}$ is an ordered triple ($V,\mathcal{G},\mathcal{B}$) such that:
\begin{adfenumerate}
\item[(i)]{$V$
is a base set of cardinality $u_1 g_1 + \dots + u_r g_r$;}
\item[(ii)]{$\mathcal{G}$
is a partition of $V$ into $u_i$ subsets of cardinality $g_i$, $i = 1, \dots, r$, called \textit{groups};}
\item[(iii)]{$\mathcal{B}$
is a non-empty collection of subsets of $V$ with cardinalities $k \in K$, called \textit{blocks}; and}
\item[(iv)]{each pair of elements from distinct groups occurs in precisely one block but no pair of
elements from the same group occurs in any block.}
\end{adfenumerate}
We often abbreviate $\{k\}$-GDD to $k$-GDD,
and a $k$-GDD of type $q^k$ is also called a {\em transversal design}, TD$(k,q)$.
A {\em parallel class} in a group divisible design is a subset of the block set in which each element of the base set appears exactly once.
A $k$-GDD is called {\em resolvable}, and is denoted by $k$-RGDD, if the entire set of blocks can be partitioned into parallel classes.
% A $k$-RGDD of type $q^k$ is also called a {\em resolvable transversal design}, RTD$(k,q)$.
If there exist $k$ mutually orthogonal Latin squares (MOLS) of side $q$, then there exists
a $(k+2)$-GDD of type $q^{k+2}$ and
a $(k+1)$-RGDD of type $q^{k+1}$, \cite[Theorem III.3.18]{AbelColbournDinitz2007}.
Furthermore, as is well known, there exist $q - 1$ MOLS of side $q$ whenever $q$ is a prime power.

Group divisible designs are widely used Design Theory,
especially in the construction of infinite classes of combinatorial designs by means of
a technique known as Wilson's Fundamental Construction,
\cite{WilsonRM1972}, \cite[Theorem IV.2.5]{GreigMullen2007}.
In general, the existence spectrum problem for group divisible designs with constant block sizes, $k$-GDDs,
appears to be a long way from being completely solved,
except of course when $k = 2$ in which case the design is essentially a complete multipartite graph.
Nevertheless, for $k = 3$ and 4, provided there are not too many distinct group sizes, considerable progress has been made.
Indeed, the spectrum has been completely determined for
3-GDDs of type $g^u$, \cite{Hanani1975}, \cite[Theorem IV.4.1]{Ge2007}, and
3-GDDs of type $g^u m^1$, \cite{ColbournHoffmanRees1992}, \cite[Theorem IV.4.4]{Ge2007}.
Furthermore, the necessary conditions for the existence of 4-GDDs of type $g^u$, namely
\begin{equation*}
\label{eqn:4-GDD g^u necessary}
u \ge 4,~~ g(u - 1) \equiv 0 \adfmod{3} \textrm{~~and~~} g^2 u (u - 1) \equiv 0 \adfmod{12},
\end{equation*}
are known to be sufficient except that there are no 4-GDDs of types $2^4$ and $6^4$, \cite{BrouwerSchrijverHanani1977}, \cite[Theorem IV.4.6]{Ge2007}.
However, in the case of 4-GDDs of type $g^u m^1$, only a partial solution has been achieved
in the sense that for each $g$ there are at most a small number of $u$ where existence remains undecided,
\cite{GeRees2002}, \cite{GeReesZhu2002}, \cite{GeLing2004}, \cite{GeRees2004}, \cite{GeLing2005}, \cite{Schuster2010},
\cite{WeiGe2013}, \cite{WeiGe2014}, \cite{Schuster2014}, \cite{WeiGe2015}, \cite{ForbesForbes2018}, \cite{Forbes2019}, \cite{Forbes2019B}.

Group divisible designs with block size 4 and type $(3a)^4 (3b)^1 (6a)^1$ have also been studied.
Rees \& Stinson \cite[Lemma 3.11]{ReesStinson1989B}
showed that the necessary conditions for the existence of these 4-GDDs, namely $1 \le a \le b \le 2a$, are sufficient,
with a few possible exceptions.
The number of possible exceptions was later reduced to four by Wang \& Shen \cite[Theorem 2.13]{WangShen2008}
and then to zero by the author of this paper \cite[Theorem 3.2]{Forbes2019}.
Also in \cite{Forbes2019} we effectively used 4-GDDs of type
$(3a)^4 (3a + 3)^1 (6a)^1$ as well as types $(6a)^6 (6a + 3)^1 (18a)^1$ and $(12a)^7 (12a + 3)^1 (42a)^1$
to complete the solution of the existence problem for 4-GDDs of type $g^u m^1$ when $g \equiv 3 \adfmod{12}$
and substantially reduce the number of possible exceptions when $g \equiv 9 \adfmod{12}$
\cite[Theorems 3.3, 3.4, 3.12]{Forbes2019}.

The main purpose of this paper is to prove that when $a$ is an integer the necessary conditions for the existence of 4-GDDs of types
$(4a)^5 b^1 (10a)^1$, $(6a+3)^6 b^1 (18a + 9)^1$ and $(6a)^6 b^1 (18a)^1$ are sufficient.
This is achieved in
Theorems~\ref{thm:4-GDD (4a)^5 b^1 (10a)^1 existence},
\ref{thm:4-GDD (6a+3)^6 b^1 (18a+9)^1 existence} and
\ref{thm:4-GDD (6a)^6 b^1 (18a)^1 existence}.

We also prove that for integer $a$, the necessary conditions suffice for 4-GDDs of type $(12a)^7 b^1 (42a)^1$,
except possibly
for $12a \in$ \{120, 180, 240, 360, 420, 720, 840\}, $24a < b < 42a$,
for $12a \in \{144, 1008\}$, $30a < b < 42a$, and
for $12a \in$ \{168, 252, 336, 504, 1512\}, $36a < b < 42a$.
This is Theorem~\ref{thm:4-GDD (12a)^7 b^1 (42a)^1 existence}.
For these values of $12a$, the range of $b$ is limited by the use of 4-GDDs of type $(6\alpha)^d \beta^1 (3\alpha d)^1$, $d \in \{4,5,6\}$
in the construction described by Theorem~\ref{thm:(6am)^d b^1 (3adm)^1}.

For the remainder of this section, we summarize some known results and constructions that we require in Section~\ref{sec:4-GDD main}.
%
%%%%%%%%%%%%%%%%%%%%%%%%%%%%%%%%%%%%%%%%%%%%%%%%%%%%%%%%%%%%%%%%%%%%%%%%%%%%%%%%%%%%%%%%%%
%%%%%%%%%%%%%%%%%%%%%%%%%%%%%%%%%%%%%%%%%%%%%%%%%%%%%%%%%%%%%%%%%%%%%%%%%%%%%%%%%%%%%%%%%%
\begin{theorem}
\label{thm:4-GDD g^u m^1 existence}
Let $g$, $u$ and $m$ be integers. There exists a 4-GDD of type $g^u m^1$ in the following cases:

{\normalfont (i)} $g \equiv 0 \adfmod{6}$, $g \ge 6$, $u \ge 4$, $m \equiv 0 \adfmod{3}$ and $0 \le m \le g(u-1)/2$
except that there is no 4-GDD of type $6^4 0^1$;

{\normalfont (ii)} $g \equiv 3 \adfmod{6}$, $g \ge 3$, $u = 7$, $m \equiv 3 \adfmod{6}$ and $0 \le m \le 3g$;

{\normalfont (iii)} $g \equiv 2$ or $4 \adfmod{6}$, $g \ge 4$, $u = 6$ and $m = 5g/2$.
\end{theorem}
{\bf Proof}~ For (i) and (ii), see \cite[Theorem 7.1]{WeiGe2015} and \cite[Theorem 1.1]{ForbesForbes2018}.
For (iii), see \cite[Theorem IV.4.11]{Ge2007}.
~\adfQED
%
%%%%%%%%%%%%%%%%%%%%%%%%%%%%%%%%%%%%%%%%%%%%%%%%%%%%%%%%%%%%%%%%%%%%%%%%%%%%%%%%%%%%%%%%%%
%%%%%%%%%%%%%%%%%%%%%%%%%%%%%%%%%%%%%%%%%%%%%%%%%%%%%%%%%%%%%%%%%%%%%%%%%%%%%%%%%%%%%%%%%%
\begin{theorem}
\label{thm:4-GDD (3a)^4 (3b)^1 (6a)^1 existence}
Let $a$ and $b$ be positive integers. There exists a 4-GDD of type $(3a)^4 (3b)^1 (6a)^1$ if and only if $1 \le a \le b \le 2a$.
\end{theorem}
{\bf Proof}~ See \cite{ReesStinson1989B}, \cite[Theorem 2.13]{WangShen2008} and \cite[Theorem 3.2]{Forbes2019}.
~\adfQED
%
%%%%%%%%%%%%%%%%%%%%%%%%%%%%%%%%%%%%%%%%%%%%%%%%%%%%%%%%%%%%%%%%%%%%%%%%%%%%%%%%%%%%%%%%%%
%%%%%%%%%%%%%%%%%%%%%%%%%%%%%%%%%%%%%%%%%%%%%%%%%%%%%%%%%%%%%%%%%%%%%%%%%%%%%%%%%%%%%%%%%%
\begin{theorem}
\label{thm:4-GDD (mg)^u}
Suppose there exists a 4-GDD of type $g_1^{u_1} \dots g_r^{u_r}$ and let $m \ge 3$, $m \neq 6$ be an integer.
Then there exists a 4-GDD of type $(mg_1)^{u_1} \dots (mg_r)^{u_r}$.
\end{theorem}
\noindent {\bf Proof~}
Take the 4-GDD of type $g_1^{u_1} \dots g_r^{u_r}$, inflate its points by a factor of $m$ and overlay the inflated blocks with 4-GDDs of type $m^4$.
~\adfQED
%
%%%%%%%%%%%%%%%%%%%%%%%%%%%%%%%%%%%%%%%%%%%%%%%%%%%%%%%%%%%%%%%%%%%%%%%%%%%%%%%%%%%%%%%%%%
%%%%%%%%%%%%%%%%%%%%%%%%%%%%%%%%%%%%%%%%%%%%%%%%%%%%%%%%%%%%%%%%%%%%%%%%%%%%%%%%%%%%%%%%%%
\begin{theorem}
\label{thm:(am)^d b^1 (adm/2)^1}
Let $\epsilon = 3$ or $6$.
Suppose there exist $d$ mutually orthogonal Latin squares of side $m$.
Suppose also there exists a 4-GDD of type $\alpha^d \beta^1 (\alpha d/2)^1$ for
each $\beta \in \{\alpha, \alpha + \epsilon, \dots, \beta_{\max}\}$, where
$\beta_{\max} \equiv \alpha \adfmod{\epsilon}$ and $\alpha \le \beta_{\max} \le \alpha d/2$.
Then there exists a 4-GDD of type $(\alpha m)^d b^1 (\alpha d m/2)^1$ for each $b \in \{\alpha m, \alpha m + \epsilon, \dots, \beta_{\max} m\}$.
\end{theorem}
{\bf Proof}~
Take a $(d + 2)$-GDD of type $m^{d + 2}$.
Choose a point, $P$, and inflate it by a factor of $\beta$ chosen from $\{\alpha, \alpha + \epsilon, \dots, \beta_{\max}\}$.
Inflate all points other than $P$ in the same group as $P$ by factors chosen from $\{\alpha, \beta_{\max}\}$.
Choose a group not containing $P$ and inflate all of its points by $\alpha d/2$.
Inflate all other points by $\alpha$.
Overlay the inflated blocks with 4-GDDs of types $\alpha^d \beta^1 (\alpha d/2)^1$, $\beta \in \{\alpha, \alpha + \epsilon, \dots, \beta_{\max}\}$.
The result is a 4-GDD of type $(\alpha m)^d b^1 (\alpha d m/2)^1$ for $b \in \{\alpha m, \alpha m + \epsilon, \dots, \beta_{\max} m\}$.
\adfQED
%
%%%%%%%%%%%%%%%%%%%%%%%%%%%%%%%%%%%%%%%%%%%%%%%%%%%%%%%%%%%%%%%%%%%%%%%%%%%%%%%%%%%%%%%%%%
%%%%%%%%%%%%%%%%%%%%%%%%%%%%%%%%%%%%%%%%%%%%%%%%%%%%%%%%%%%%%%%%%%%%%%%%%%%%%%%%%%%%%%%%%%
\begin{theorem}
\label{thm:(6am)^d b^1 (3adm)^1}
Let $\epsilon = 3$ or $6$. % $\epsilon = 6$ might not be necessary
Let $m$, $d$ and $\alpha$ be positive integers with $5 \le m \le d + 1$.
Suppose there exist $m - 1$ mutually orthogonal Latin squares of side $d + 1$.
Suppose also there exists a 4-GDD of type $(6\alpha)^{m - 1} \beta^1 (3\alpha(m-1))^1$ and $(6\alpha)^d \beta^1 (3\alpha d)^1$
for each $\beta \in \{6\alpha, 6\alpha + \epsilon, \dots, \beta_{\max}\}$,
where $\beta_{\max} \equiv 0 \adfmod{\epsilon}$ and $6\alpha \le \beta_{\max} \le 3\alpha(m - 1)$.
Then there exists a 4-GDD of type $(6\alpha m)^d b^1 (3\alpha m d)^1$ for each $b \equiv 0 \adfmod{\epsilon}$ and $6\alpha \le b \le \beta_{\max}m$.
\end{theorem}
{\bf Proof}~ Take an $(m + 1)$-GDD of type $(d + 1)^{m + 1}$,
choose a point, $P$, and remove the blocks of the parallel class defined by $P$.
Choose one of the removed blocks and inflate its points other than $P$ by factors
$\beta_1$, $\beta_2$, \dots, $\beta_m$, each chosen from $\{6\alpha, 6\alpha + \epsilon, \dots, \beta_{\max}\}$.
Inflate all points other than $P$ in the same group as $P$ by a factor of $3\alpha (m - 1)$.
Inflate $P$ by a factor of $3\alpha d$, and
inflate the remaining points by $6\alpha$.

Overlay the surviving inflated blocks with 4-GDDs of type $(6\alpha)^{m - 1} \beta^1 (3\alpha(m - 1))^1$.
% Include $(6\alpha)^m (3\alpha(m - 1))^1$ in case $m \le d$ and there is no $\beta = 6 \alpha$.
Take each inflated group other than the one containing $P$ and overlay it plus the inflated point $P$ with a
4-GDD of type $(6\alpha)^d \beta^1 (3\alpha d)^1$.
In each case $\beta \in \{6\alpha, 6\alpha + \epsilon, \dots, \beta_{\max}\}$.
The result is a 4-GDD of type $(6\alpha m)^d b^1 (3\alpha m d)^1$, where
$b = \sum_{i=1}^m \beta_i$ can take any value in the set $\{6\alpha m, 6\alpha m + \epsilon, \dots, \beta_{\max}m\}$.
\adfQED

\adfPgap
For the next theorem, we require the concept of a transversal design with a hole, which is defined as follows.
An {\em incomplete transversal design}, ITD$(k, g; h)$, with block size $k$, group size $g$ and hole size $h k$,
is a quadruple $(V, \mathcal{G}, H, \mathcal{B})$, where:
$V$ is a set of $kg$ elements;
$\mathcal{G}$ is a partition of $V$ into $k$ groups, each of size $g$;
$H$, the hole, is a subset of $V$ such that $|G \cap H| = h$ for each $G \in \mathcal{G}$; and
$\mathcal{B}$ is a collection of $k$-subsets of $V$, the blocks, with the property that
(i) every unordered pair of elements of $V$ that does not occur in a group or in $H$ is contained in precisely one block, and
(ii) no pair of elements of a group or of $H$ occurs in any block.

For a discussion of the related concept of incomplete mutually orthogonal Latin squares, IMOLS,
we refer the reader to \cite{AbelColbournDinitz2007b}.
Here, is convenient to recall that for $k \ge 4$, an ITD$(k, g; h)$ exists if and only if there exist
$k - 2$ IMOLS of side $g$ and hole size $h \times h$, \cite[Theorem III.4.11]{AbelColbournDinitz2007b},
and for the existence of these structures we cite \cite[Table III.4.14]{AbelColbournDinitz2007b}.
%
%%%%%%%%%%%%%%%%%%%%%%%%%%%%%%%%%%%%%%%%%%%%%%%%%%%%%%%%%%%%%%%%%%%%%%%%%%%%%%%%%%%%%%%%%%
%%%%%%%%%%%%%%%%%%%%%%%%%%%%%%%%%%%%%%%%%%%%%%%%%%%%%%%%%%%%%%%%%%%%%%%%%%%%%%%%%%%%%%%%%%
\begin{theorem}
\label{thm:(am)^d b^1 (adm/2)^1 hole}
Let $\epsilon = 3$ or $6$.
Suppose there exists an ITD$(d + 2, m; h)$.
Suppose also there exists a 4-GDD of type $\alpha^d \beta^1 (\alpha d/2)^1$ for
each $\beta \in \{\alpha, \alpha + \epsilon, \dots, \beta_{\max}\}$, where
$\beta_{\max} \equiv \alpha \adfmod{\epsilon}$ and $\alpha \le \beta_{\max} \le \alpha d/2$.
Suppose furthermore there exists a 4-GDD of type $(\alpha h)^d \gamma^1 (\alpha h d/2)^1$ for
$\gamma \in \{\alpha h, \alpha h + \epsilon, \dots, \beta_{\max} h\}$.
Then there exists a 4-GDD of type $(\alpha m)^d b^1 (\alpha d m/2)^1$ for each $b \in \{\alpha m, \alpha m + \epsilon, \dots, \beta_{\max} m\}$.
\end{theorem}
{\bf Proof}~
Take an ITD$(d + 2, m; h)$.
Choose a group, $G$, and inflate its points by factors $\beta$ chosen from $\{\alpha, \alpha + \epsilon, \dots, \beta_{\max}\}$.
Choose another group and inflate each of its points by $\alpha d/2$.
Inflate all other points by $\alpha$.

Overlay the inflated blocks with 4-GDDs of types $\alpha^d \beta^1 (\alpha d/2)^1$, $\beta \in \{\alpha, \alpha + \epsilon, \dots, \beta_{\max}\}$.
Overlay the hole with a 4-GDD of type $\alpha h^d \gamma^1 (\alpha h d/2)^1$, where
$\gamma$ is the sum of the weights $\beta$ over the $h$ points that occur in the intersection of $G$ and the hole.
The result is a 4-GDD of type $(\alpha m)^d b^1 (\alpha d m/2)^1$ for $b \in \{\alpha m, \alpha m + \epsilon, \dots, \beta_{\max} m\}$.
\adfQED

%%%%%%%%%%%%%%%%%%%%%%%%%%%%%%%%%%%%%%%%%%%%%%%%%%%%%%%%%%%%%%%%%%%%%%%%%%%%%%%%%%%%%%%%%%
%%%%%%%%%%%%%%%%%%%%%%%%%%%%%%%%%%%%%%%%%%%%%%%%%%%%%%%%%%%%%%%%%%%%%%%%%%%%%%%%%%%%%%%%%%
%%%%%%%%%%%%%%%%%%%%%%%%%%%%%%%%%%%%%%%%%%%%%%%%%%%%%%%%%%%%%%%%%%%%%%%%%%%%%%%%%%%%%%%%%%
%%%%%%%%%%%%%%%%%%%%%%%%%%%%%%%%%%%%%%%%%%%%%%%%%%%%%%%%%%%%%%%%%%%%%%%%%%%%%%%%%%%%%%%%%%

\section{New 4-GDDs}\label{sec:new 4-GDDs}

The 4-GDDs that we use to prove our theorems are given here, grouped into lemmas.
The blocks of the designs are generated by appropriate mappings from sets of base blocks,
which are collected in the appendix together with the instructions for expanding them.
The presentation has not been refined in any way; the base blocks are listed exactly as created by the author's computer system.
It is possible that some of the designs might be constructible by other means.
If the point set of the 4-GDD has $v$ elements, it is represented by $Z_v = \{0, 1, \dots, v-1\}$ partitioned into groups as indicated.
The expression $a\adfmod{b}$ denotes the integer $n$ such that $0 \le n < b$ and $b \,|\, n - a$.

As an aid to checking the correctness of the designs we also provide the block generation instructions in a compact format: $(v,((s,b/s,M)),T)$,
where $v$ is the number of points, $s$ is the number of base blocks, $b$ is the number of blocks,
$M$ encodes the mapping and $T$ encodes the group type.
We leave it for the interested reader to determine the structure of $M$ and $T$ from the corresponding English descriptions.
% In case of inconsistency please ignore the English form.
This feature might be especially useful for the larger 4-GDDs,
those of Lemma~\ref{lem:4-GDD 180^6 b^1 540^1}, for example, where each design is generated from over 1000 base blocks.

\newcommand{\adfAppRefA}[2]{#1}% Reference to Appendix section in the arXiv paper.

\newcommand{\adfAppRefJ}[2]{#2}% Reference to Appendix section in the journal paper.

%%%%%%%%%%%%%%%%%%%%%%%%%%%%%%%%%%%%%%%%%%%%%%%
%%%%%%%%%%%%%%%%%%%%%%%%%%%%%%%%%%%%%%%%%%%%%%%
\begin{lemma}
\label{lem:4-GDD 4^5 b^1 10^1}
\label{lem:4-GDD 8^5 b^1 20^1}
\label{lem:4-GDD 10^5 b^1 25^1}
\label{lem:4-GDD 12^5 b^1 30^1}
There exist 4-GDDs of types
$4^5  7^1 10^1$,
$8^5  11^1 20^1$,
$8^5  17^1 20^1$,
$8^5  20^2$,
$10^5 16^1 25^1$,
$10^5 22^1 25^1$,
$12^5 18^1 30^1$,
$12^5 24^1 30^1$,
$12^5 27^1 30^1$.
\end{lemma}
{\bf Proof}~ The designs are presented in Appendix~\adfAppRefJ{\ref{app:4-GDD 4^5 b^1 10^1}}{A}.
~\adfQED

%%%%%%%%%%%%%%%%%%%%%%%%%%%%%%%%%%%%%%%%%%%%%%%
%%%%%%%%%%%%%%%%%%%%%%%%%%%%%%%%%%%%%%%%%%%%%%%
\begin{lemma}
\label{lem:4-GDD 16^5 b^1 40^1}
\label{lem:4-GDD 20^5 b^1 50^1}
There exist 4-GDDs of types
$16^5 19^1 40^1$,
$16^5 22^1 40^1$,
$16^5 25^1 40^1$,
$16^5 31^1 40^1$,
$16^5 34^1 40^1$,
$16^5 37^1 40^1$,
$20^5 23^1 50^1$,
$20^5 26^1 50^1$,
$20^5 29^1 50^1$,
$20^5 32^1 50^1$,
$20^5 38^1 50^1$,
$20^5 41^1 50^1$,
$20^5 44^1 50^1$,
$20^5 47^1 50^1$.
\end{lemma}
{\bf Proof}~ The designs are presented in Appendix~\adfAppRefA{\ref{app:4-GDD 16^5 b^1 40^1}}{B}.
~\adfQED

%%%%%%%%%%%%%%%%%%%%%%%%%%%%%%%%%%%%%%%%%%%%%%%
%%%%%%%%%%%%%%%%%%%%%%%%%%%%%%%%%%%%%%%%%%%%%%%
\begin{lemma}
\label{lem:4-GDD 24^5 b^1 60^1}
\label{lem:4-GDD 30^5 b^1 75^1}
There exist 4-GDDs of types
$24^5 27^1 60^1$,
$24^5 30^1 60^1$,
$24^5 36^1 60^1$,
$24^5 39^1 60^1$,
$24^5 45^1 60^1$,
$24^5 54^1 60^1$,
$24^5 57^1 60^1$,
$30^5 36^1 75^1$,
$30^5 42^1 75^1$,
$30^5 54^1 75^1$,
$30^5 72^1 75^1$.
\end{lemma}
{\bf Proof}~ The designs are presented in Appendix~\adfAppRefA{\ref{app:4-GDD 24^5 b^1 60^1}}{C}.
~\adfQED

%%%%%%%%%%%%%%%%%%%%%%%%%%%%%%%%%%%%%%%%%%%%%%%
%%%%%%%%%%%%%%%%%%%%%%%%%%%%%%%%%%%%%%%%%%%%%%%
\begin{lemma}
\label{lem:4-GDD 40^5 b^1 100^1}
There exist 4-GDDs of types
$40^5 43^1 100^1$,
$40^5 46^1 100^1$,
$40^5 49^1 100^1$,
$40^5 52^1 100^1$,
$40^5 58^1 100^1$,
$40^5 61^1 100^1$,
$40^5 67^1 100^1$,
$40^5 73^1 100^1$,
$40^5 76^1 100^1$,
$40^5 79^1 100^1$,
$40^5 82^1 100^1$,
$40^5 91^1 100^1$,
$40^5 94^1 100^1$,
$40^5 97^1 100^1$.
\end{lemma}
{\bf Proof}~ The designs are presented in Appendix~\adfAppRefA{\ref{app:4-GDD 40^5 b^1 100^1}}{D}.
~\adfQED

%%%%%%%%%%%%%%%%%%%%%%%%%%%%%%%%%%%%%%%%%%%%%%%
%%%%%%%%%%%%%%%%%%%%%%%%%%%%%%%%%%%%%%%%%%%%%%%
\begin{lemma}
\label{lem:4-GDD 60^5 b^1 150^1}
There exist 4-GDDs of types
$60^5 63^1 150^1$,
$60^5 66^1 150^1$,
$60^5 72^1 150^1$,
$60^5 81^1 150^1$,
$60^5 84^1 150^1$,
$60^5 93^1 150^1$,
$60^5 99^1 150^1$,
$60^5 102^1 150^1$,
$60^5 108^1 150^1$,
$60^5 111^1 150^1$,
$60^5 117^1 150^1$,
$60^5 126^1 150^1$,
$60^5 129^1 150^1$,
$60^5 138^1 150^1$,
$60^5 144^1 150^1$,
$60^5 147^1 150^1$.
\end{lemma}
{\bf Proof}~ The designs are presented in Appendix~\adfAppRefA{\ref{app:4-GDD 60^5 b^1 150^1}}{E}.
~\adfQED

%%%%%%%%%%%%%%%%%%%%%%%%%%%%%%%%%%%%%%%%%%%%%%%
%%%%%%%%%%%%%%%%%%%%%%%%%%%%%%%%%%%%%%%%%%%%%%%
\begin{lemma}
\label{lem:4-GDD 80^5 b^1 200^1}
There exist 4-GDDs of types
$80^5 83^1 200^1$,
$80^5 89^1 200^1$,
$80^5 101^1 200^1$,
$80^5 107^1 200^1$,
$80^5 113^1 200^1$,
$80^5 119^1 200^1$,
$80^5 131^1 200^1$,
$80^5 137^1 200^1$,
$80^5 143^1 200^1$,
$80^5 149^1 200^1$,
$80^5 161^1 200^1$,
$80^5 167^1 200^1$,
$80^5 173^1 200^1$,
$80^5 179^1 200^1$,
$80^5 182^1 200^1$,
$80^5 191^1 200^1$,
$80^5 194^1 200^1$,
$80^5 197^1 200^1$.
\end{lemma}
{\bf Proof}~ The designs are presented in Appendix~\adfAppRefA{\ref{app:4-GDD 80^5 b^1 200^1}}{F}.
~\adfQED

%%%%%%%%%%%%%%%%%%%%%%%%%%%%%%%%%%%%%%%%%%%%%%%
%%%%%%%%%%%%%%%%%%%%%%%%%%%%%%%%%%%%%%%%%%%%%%%
\begin{lemma}
\label{lem:4-GDD 6^6 b^1 18^1}
\label{lem:4-GDD 9^6 b^1 27^1}
\label{lem:4-GDD 12^6 b^1 36^1}
There exist 4-GDDs of types
$6^6 12^1 18^1$,
$6^6 15^1 18^1$,
$6^6 18^2$,
$9^6 21^1 27^1$,
$12^6 18^1 36^1$,
$12^6 24^1 36^1$,
$12^6 27^1 36^1$,
$12^6 30^1 36^1$,
$12^6 33^1 36^1$,
$12^6 36^2$.
\end{lemma}
{\bf Proof}~ The designs are presented in Appendix~\adfAppRefA{\ref{app:4-GDD 6^6 b^1 18^1}}{G}.
~\adfQED

%%%%%%%%%%%%%%%%%%%%%%%%%%%%%%%%%%%%%%%%%%%%%%%
%%%%%%%%%%%%%%%%%%%%%%%%%%%%%%%%%%%%%%%%%%%%%%%
\begin{lemma}
\label{lem:4-GDD 15^6 b^1 45^1}
\label{lem:4-GDD 18^6 b^1 54^1}
There exist 4-GDDs of types
$15^6 21^1 45^1$,
$15^6 27^1 45^1$,
$15^6 33^1 45^1$,
$15^6 39^1 45^1$,
$18^6 24^1 54^1$,
$18^6 30^1 54^1$,
$18^6 33^1 54^1$,
$18^6 39^1 54^1$,
$18^6 42^1 54^1$,
$18^6 48^1 54^1$,
$18^6 51^1 54^1$.
\end{lemma}
{\bf Proof}~ The designs are presented in Appendix~\adfAppRefA{\ref{app:4-GDD 15^6 b^1 45^1}}{H}
~\adfQED

%%%%%%%%%%%%%%%%%%%%%%%%%%%%%%%%%%%%%%%%%%%%%%%
%%%%%%%%%%%%%%%%%%%%%%%%%%%%%%%%%%%%%%%%%%%%%%%
\begin{lemma}
\label{lem:4-GDD 24^6 b^1 72^1}
There exist 4-GDDs of types
$24^6 30^1 72^1$,
$24^6 33^1 72^1$,
$24^6 39^1 72^1$,
$24^6 42^1 72^1$,
$24^6 45^1 72^1$,
$24^6 51^1 72^1$,
$24^6 54^1 72^1$,
$24^6 57^1 72^1$,
$24^6 63^1 72^1$,
$24^6 66^1 72^1$,
$24^6 69^1 72^1$.
\end{lemma}
{\bf Proof}~ The designs are presented in Appendix~\adfAppRefA{\ref{app:4-GDD 24^6 b^1 72^1}}{I}.
~\adfQED

%%%%%%%%%%%%%%%%%%%%%%%%%%%%%%%%%%%%%%%%%%%%%%%
%%%%%%%%%%%%%%%%%%%%%%%%%%%%%%%%%%%%%%%%%%%%%%%
\begin{lemma}
\label{lem:4-GDD 30^6 b^1 90^1}
There exist 4-GDDs of types
$30^6 63^1 90^1$,
$30^6 66^1 90^1$,
$30^6 69^1 90^1$,
$30^6 72^1 90^1$,
$30^6 78^1 90^1$,
$30^6 81^1 90^1$,
$30^6 84^1 90^1$,
$30^6 87^1 90^1$.
\end{lemma}
{\bf Proof}~ The designs are presented in Appendix~\adfAppRefA{\ref{app:4-GDD 30^6 b^1 90^1}}{J}.
~\adfQED

%%%%%%%%%%%%%%%%%%%%%%%%%%%%%%%%%%%%%%%%%%%%%%%
%%%%%%%%%%%%%%%%%%%%%%%%%%%%%%%%%%%%%%%%%%%%%%%
\begin{lemma}
\label{lem:4-GDD 36^6 b^1 108^1}
There exist 4-GDDs of types
$36^6 39^1 108^1$,
$36^6 42^1 108^1$,
$36^6 48^1 108^1$,
$36^6 51^1 108^1$,
$36^6 57^1 108^1$,
$36^6 60^1 108^1$,
$36^6 66^1 108^1$,
$36^6 69^1 108^1$,
$36^6 75^1 108^1$,
$36^6 78^1 108^1$,
$36^6 84^1 108^1$,
$36^6 87^1 108^1$,
$36^6 93^1 108^1$,
$36^6 96^1 108^1$,
$36^6 102^1 108^1$,
$36^6 105^1 108^1$.
\end{lemma}
{\bf Proof}~ The designs are presented in Appendix~\adfAppRefA{\ref{app:4-GDD 36^6 b^1 108^1}}{K}.
~\adfQED

%%%%%%%%%%%%%%%%%%%%%%%%%%%%%%%%%%%%%%%%%%%%%%%
%%%%%%%%%%%%%%%%%%%%%%%%%%%%%%%%%%%%%%%%%%%%%%%
\begin{lemma}
\label{lem:4-GDD 45^6 b^1 135^1}
There exist 4-GDDs of types
$45^6 51^1 135^1$,
$45^6 57^1 135^1$,
$45^6 69^1 135^1$,
$45^6 87^1 135^1$,
$45^6 93^1 135^1$,
$45^6 111^1 135^1$,
$45^6 123^1 135^1$,
$45^6 129^1 135^1$.
\end{lemma}
{\bf Proof}~ The designs are presented in Appendix~\adfAppRefA{\ref{app:4-GDD 45^6 b^1 135^1}}{L}.
~\adfQED

%%%%%%%%%%%%%%%%%%%%%%%%%%%%%%%%%%%%%%%%%%%%%%%
%%%%%%%%%%%%%%%%%%%%%%%%%%%%%%%%%%%%%%%%%%%%%%%
\begin{lemma}
\label{lem:4-GDD 60^6 b^1 180^1}
There exist 4-GDDs of types
$60^6 123^1 180^1$,
$60^6 126^1 180^1$,
$60^6 129^1 180^1$,
$60^6 138^1 180^1$,
$60^6 141^1 180^1$,
$60^6 144^1 180^1$,
$60^6 147^1 180^1$,
$60^6 153^1 180^1$,
$60^6 159^1 180^1$,
$60^6 162^1 180^1$,
$60^6 168^1 180^1$,
$60^6 171^1 180^1$,
$60^6 174^1 180^1$,
$60^6 177^1 180^1$.
\end{lemma}
{\bf Proof}~ The designs are presented in Appendix~\adfAppRefA{\ref{app:4-GDD 60^6 b^1 180^1}}{M}.
~\adfQED

%%%%%%%%%%%%%%%%%%%%%%%%%%%%%%%%%%%%%%%%%%%%%%%
%%%%%%%%%%%%%%%%%%%%%%%%%%%%%%%%%%%%%%%%%%%%%%%
\begin{lemma}
\label{lem:4-GDD 72^6 b^1 216^1}
There exist 4-GDDs of types
$72^6 183^1 216^1$,
$72^6 195^1 216^1$,
$72^6 201^1 216^1$,
$72^6 213^1 216^1$.
\end{lemma}
{\bf Proof}~ The designs are presented in Appendix~\adfAppRefA{\ref{app:4-GDD 72^6 b^1 216^1}}{N}.
~\adfQED

%%%%%%%%%%%%%%%%%%%%%%%%%%%%%%%%%%%%%%%%%%%%%%%
%%%%%%%%%%%%%%%%%%%%%%%%%%%%%%%%%%%%%%%%%%%%%%%
\begin{lemma}
\label{lem:4-GDD 90^6 b^1 270^1}
There exist 4-GDDs of types
$90^6 183^1 270^1$,
$90^6 186^1 270^1$,
$90^6 192^1 270^1$,
$90^6 201^1 270^1$,
$90^6 204^1 270^1$,
$90^6 213^1 270^1$,
$90^6 219^1 270^1$,
$90^6 222^1 270^1$,
$90^6 228^1 270^1$,
$90^6 231^1 270^1$,
$90^6 237^1 270^1$,
$90^6 246^1 270^1$,
$90^6 249^1 270^1$,
$90^6 258^1 270^1$,
$90^6 264^1 270^1$,
$90^6 267^1 270^1$.
\end{lemma}
{\bf Proof}~ The designs are presented in Appendix~\adfAppRefA{\ref{app:4-GDD 90^6 b^1 270^1}}{O}.
~\adfQED

%%%%%%%%%%%%%%%%%%%%%%%%%%%%%%%%%%%%%%%%%%%%%%%
%%%%%%%%%%%%%%%%%%%%%%%%%%%%%%%%%%%%%%%%%%%%%%%
\begin{lemma}
\label{lem:4-GDD 120^6 b^1 360^1}
There exist 4-GDDs of types
$120^6 243^1 360^1$,
$120^6 249^1 360^1$,
$120^6 261^1 360^1$,
$120^6 267^1 360^1$,
$120^6 273^1 360^1$,
$120^6 279^1 360^1$,
$120^6 291^1 360^1$,
$120^6 297^1 360^1$,
$120^6 303^1 360^1$,
$120^6 309^1 360^1$,
$120^6 321^1 360^1$,
$120^6 327^1 360^1$,
$120^6 333^1 360^1$,
$120^6 339^1 360^1$,
$120^6 351^1 360^1$,
$120^6 357^1 360^1$.
\end{lemma}
{\bf Proof}~ The designs are presented in Appendix~\adfAppRefA{\ref{app:4-GDD 120^6 b^1 360^1}}{P}.
~\adfQED

%%%%%%%%%%%%%%%%%%%%%%%%%%%%%%%%%%%%%%%%%%%%%%%
%%%%%%%%%%%%%%%%%%%%%%%%%%%%%%%%%%%%%%%%%%%%%%%
\begin{lemma}
\label{lem:4-GDD 180^6 b^1 540^1}
There exist 4-GDDs of types
$180^6 363^1 540^1$,
$180^6 381^1 540^1$,
$180^6 393^1 540^1$,
$180^6 399^1 540^1$,
$180^6 411^1 540^1$,
$180^6 417^1 540^1$,
$180^6 429^1 540^1$,
$180^6 447^1 540^1$,
$180^6 453^1 540^1$,
$180^6 471^1 540^1$,
$180^6 483^1 540^1$,
$180^6 489^1 540^1$,
$180^6 501^1 540^1$,
$180^6 507^1 540^1$,
$180^6 519^1 540^1$,
$180^6 537^1 540^1$.
\end{lemma}
{\bf Proof}~ The designs are presented in Appendix~\adfAppRefA{\ref{app:4-GDD 180^6 b^1 540^1}}{Q}.
~\adfQED

%%%%%%%%%%%%%%%%%%%%%%%%%%%%%%%%%%%%%%%%%%%%%%%
%%%%%%%%%%%%%%%%%%%%%%%%%%%%%%%%%%%%%%%%%%%%%%%
\begin{lemma}
\label{lem:4-GDD 6^7 b^1 21^1}
\label{lem:4-GDD 12^7 b^1 42^1}
\label{lem:4-GDD 18^7 b^1 63^1}
There exist 4-GDDs of types
$6^7 12^1 21^1$,
$6^7 18^1 21^1$,
$12^7 18^1 42^1$,
$12^7 21^1 42^1$,
$12^7 24^1 42^1$,
$12^7 27^1 42^1$,
$12^7 30^1 42^1$,
$12^7 33^1 42^1$,
$12^7 36^1 42^1$,
$12^7 39^1 42^1$,
$12^7 42^2$,
$18^7 48^1 63^1$,
$18^7 60^1 63^1$.
\end{lemma}
{\bf Proof}~ The designs are presented in Appendix~\adfAppRefA{\ref{app:4-GDD 6^7 b^1 21^1}}{R}.
~\adfQED

%%%%%%%%%%%%%%%%%%%%%%%%%%%%%%%%%%%%%%%%%%%%%%%
%%%%%%%%%%%%%%%%%%%%%%%%%%%%%%%%%%%%%%%%%%%%%%%
\begin{lemma}
\label{lem:4-GDD 24^7 b^1 84^1}
There exist 4-GDDs of types
$24^7 30^1 84^1$,
$24^7 33^1 84^1$,
$24^7 36^1 84^1$,
$24^7 39^1 84^1$,
$24^7 42^1 84^1$,
$24^7 45^1 84^1$,
$24^7 51^1 84^1$,
$24^7 54^1 84^1$,
$24^7 57^1 84^1$,
$24^7 60^1 84^1$,
$24^7 63^1 84^1$,
$24^7 66^1 84^1$,
$24^7 69^1 84^1$,
$24^7 75^1 84^1$,
$24^7 78^1 84^1$,
$24^7 81^1 84^1$,
$24^7 84^2$.
\end{lemma}
{\bf Proof}~ The designs are presented in Appendix~\adfAppRefA{\ref{app:4-GDD 24^7 b^1 84^1}}{S}.
~\adfQED

%%%%%%%%%%%%%%%%%%%%%%%%%%%%%%%%%%%%%%%%%%%%%%%
%%%%%%%%%%%%%%%%%%%%%%%%%%%%%%%%%%%%%%%%%%%%%%%
\begin{lemma}
\label{lem:4-GDD 36^7 b^1 126^1}
There exist 4-GDDs of types
$36^7 42^1 126^1$,
$36^7 48^1 126^1$,
$36^7 51^1 126^1$,
$36^7 57^1 126^1$,
$36^7 60^1 126^1$,
$36^7 66^1 126^1$,
$36^7 69^1 126^1$,
$36^7 75^1 126^1$,
$36^7 78^1 126^1$,
$36^7 84^1 126^1$,
$36^7 87^1 126^1$,
$36^7 93^1 126^1$,
$36^7 96^1 126^1$,
$36^7 102^1 126^1$,
$36^7 105^1 126^1$,
$36^7 111^1 126^1$,
$36^7 114^1 126^1$,
$36^7 120^1 126^1$,
$36^7 123^1 126^1$.
\end{lemma}
{\bf Proof}~ The designs are presented in Appendix~\adfAppRefA{\ref{app:4-GDD 36^7 b^1 126^1}}{T}.
~\adfQED

%%%%%%%%%%%%%%%%%%%%%%%%%%%%%%%%%%%%%%%%%%%%%%%
%%%%%%%%%%%%%%%%%%%%%%%%%%%%%%%%%%%%%%%%%%%%%%%
\begin{lemma}
\label{lem:4-GDD 48^7 b^1 168^1}
There exist 4-GDDs of types
$48^7 54^1 168^1$,
$48^7 57^1 168^1$,
$48^7 63^1 168^1$,
$48^7 66^1 168^1$,
$48^7 69^1 168^1$,
$48^7 75^1 168^1$,
$48^7 78^1 168^1$,
$48^7 81^1 168^1$,
$48^7 87^1 168^1$,
$48^7 90^1 168^1$,
$48^7 93^1 168^1$,
$48^7 99^1 168^1$,
$48^7 102^1 168^1$,
$48^7 105^1 168^1$,
$48^7 111^1 168^1$,
$48^7 114^1 168^1$,
$48^7 117^1 168^1$,
$48^7 123^1 168^1$,
$48^7 126^1 168^1$,
$48^7 129^1 168^1$,
$48^7 135^1 168^1$,
$48^7 138^1 168^1$,
$48^7 141^1 168^1$,
$48^7 147^1 168^1$,
$48^7 150^1 168^1$,
$48^7 153^1 168^1$,
$48^7 159^1 168^1$,
$48^7 162^1 168^1$,
$48^7 165^1 168^1$.
\end{lemma}
{\bf Proof}~ The designs are presented in Appendix~\adfAppRefA{\ref{app:4-GDD 48^7 b^1 168^1}}{U}.
~\adfQED

%%%%%%%%%%%%%%%%%%%%%%%%%%%%%%%%%%%%%%%%%%%%%%%
%%%%%%%%%%%%%%%%%%%%%%%%%%%%%%%%%%%%%%%%%%%%%%%
\begin{lemma}
\label{lem:4-GDD 60^7 b^1 210^1}
There exist 4-GDDs of types
$60^7 123^1 210^1$,
$60^7 126^1 210^1$,
$60^7 129^1 210^1$,
$60^7 132^1 210^1$,
$60^7 138^1 210^1$,
$60^7 141^1 210^1$,
$60^7 144^1 210^1$,
$60^7 147^1 210^1$,
$60^7 153^1 210^1$,
$60^7 156^1 210^1$,
$60^7 159^1 210^1$,
$60^7 162^1 210^1$,
$60^7 168^1 210^1$,
$60^7 171^1 210^1$,
$60^7 174^1 210^1$,
$60^7 177^1 210^1$,
$60^7 183^1 210^1$,
$60^7 186^1 210^1$,
$60^7 189^1 210^1$,
$60^7 192^1 210^1$,
$60^7 198^1 210^1$,
$60^7 201^1 210^1$,
$60^7 204^1 210^1$,
$60^7 207^1 210^1$.
\end{lemma}
{\bf Proof}~ The designs are presented in Appendix~\adfAppRefA{\ref{app:4-GDD 60^7 b^1 210^1}}{V}.
~\adfQED

%%%%%%%%%%%%%%%%%%%%%%%%%%%%%%%%%%%%%%%%%%%%%%%
%%%%%%%%%%%%%%%%%%%%%%%%%%%%%%%%%%%%%%%%%%%%%%%
\begin{lemma}
\label{lem:4-GDD 72^7 b^1 252^1}
There exist 4-GDDs of types
$72^7 183^1 252^1$,
$72^7 186^1 252^1$,
$72^7 195^1 252^1$,
$72^7 201^1 252^1$,
$72^7 204^1 252^1$,
$72^7 210^1 252^1$,
$72^7 213^1 252^1$,
$72^7 219^1 252^1$,
$72^7 222^1 252^1$,
$72^7 228^1 252^1$,
$72^7 231^1 252^1$,
$72^7 237^1 252^1$,
$72^7 246^1 252^1$,
$72^7 249^1 252^1$.
\end{lemma}
{\bf Proof}~ The designs are presented in Appendix~\adfAppRefA{\ref{app:4-GDD 72^7 b^1 252^1}}{W}.
~\adfQED

%%%%%%%%%%%%%%%%%%%%%%%%%%%%%%%%%%%%%%%%%%%%%%%
%%%%%%%%%%%%%%%%%%%%%%%%%%%%%%%%%%%%%%%%%%%%%%%
\begin{lemma}
\label{lem:4-GDD 84^7 b^1 294^1}
There exist 4-GDDs of types
$84^7 255^1 294^1$,
$84^7 258^1 294^1$,
$84^7 261^1 294^1$,
$84^7 264^1 294^1$,
$84^7 267^1 294^1$,
$84^7 270^1 294^1$,
$84^7 276^1 294^1$,
$84^7 279^1 294^1$,
$84^7 282^1 294^1$,
$84^7 285^1 294^1$,
$84^7 288^1 294^1$,
$84^7 291^1 294^1$.
\end{lemma}
{\bf Proof}~ The designs are presented in Appendix~\adfAppRefA{\ref{app:4-GDD 84^7 b^1 294^1}}{X}.
~\adfQED

%%%%%%%%%%%%%%%%%%%%%%%%%%%%%%%%%%%%%%%%%%%%%%%%%%%%%%%%%%%%%%%%%%%%%%%%%%%%%%%%%%%%%%%%%%
%%%%%%%%%%%%%%%%%%%%%%%%%%%%%%%%%%%%%%%%%%%%%%%%%%%%%%%%%%%%%%%%%%%%%%%%%%%%%%%%%%%%%%%%%%
%%%%%%%%%%%%%%%%%%%%%%%%%%%%%%%%%%%%%%%%%%%%%%%%%%%%%%%%%%%%%%%%%%%%%%%%%%%%%%%%%%%%%%%%%%
%%%%%%%%%%%%%%%%%%%%%%%%%%%%%%%%%%%%%%%%%%%%%%%%%%%%%%%%%%%%%%%%%%%%%%%%%%%%%%%%%%%%%%%%%%

\section{Main results}
\label{sec:4-GDD main}
%
%%%%%%%%%%%%%%%%%%%%%%%%%%%%%%%%%%%%%%%%%%%%%%%%%%%%%%%%%%%%%%%%%%%%%%%%%%%%%%%%%%%%%%%%%%
%%%%%%%%%%%%%%%%%%%%%%%%%%%%%%%%%%%%%%%%%%%%%%%%%%%%%%%%%%%%%%%%%%%%%%%%%%%%%%%%%%%%%%%%%%
\begin{theorem}
\label{thm:(4a)^5 b^1 (10a)^1 small}
There exists a 4-GDD of type $g^5 b^1 (5g/2)^1$ for
$g \in$ \{$4$, $8$, $12$, $16$, $20$, $24$, $40$, $60$, $80$, $120$\} and $b \in \{g, g + 3, \dots, 5g/2\}$.

There exists a 4-GDD of type $g^5 b^1 (5g/2)^1$ for
$g \in$ \{$10$, $30$\} and $b \in \{g, g + 6, \dots, 5g/2 - 3\}$.
\end{theorem}
{\bf Proof}~ By Theorem~\ref{thm:4-GDD g^u m^1 existence}, we may assume $b > g$.
We deal with each value of $g$ in turn.

{\boldmath $g = 4$}~
For $b = 10$, see \cite[Lemma 2.1]{Forbes2019B}.

For $b = 7$, see Lemma~\ref{lem:4-GDD 4^5 b^1 10^1}.

{\boldmath $g = 8$}~
For $b = 14$, see \cite[Lemma 2.1]{Forbes2019B}.

For $b \in \{11, 17, 20\}$, see Lemma~\ref{lem:4-GDD 8^5 b^1 20^1}.

{\boldmath $g = 10$}~
For $b \in \{16, 22\}$, see Lemma~\ref{lem:4-GDD 10^5 b^1 25^1}.

{\boldmath $g = 12$}~ For $b \in \{15, 21\}$, see \cite[Proposition IV.4.15]{Ge2007}.

For $b = 30$, use Theorem~\ref{thm:4-GDD (mg)^u} with a 4-GDD of type $4^5 10^1 10^1$.

For $b \in \{18, 24, 27\}$, see Lemma~\ref{lem:4-GDD 12^5 b^1 30^1}.

{\boldmath $g = 16$}~
For $b \in \{28, 40\}$, use Theorem~\ref{thm:4-GDD (mg)^u} with a 4-GDD of type $4^5 (b/4)^1 10^1$.

For $b \in \{19, 22, 25, 31, 34, 37\}$, see Lemma~\ref{lem:4-GDD 16^5 b^1 40^1}.

{\boldmath $g = 20$}~
For $b \in \{35, 50\}$, use Theorem~\ref{thm:4-GDD (mg)^u} with a 4-GDD of type $4^5 (b/5)^1 10^1$.

For $b \in \{23, 26, 29, 32, 38, 41, 44, 47\}$, see Lemma~\ref{lem:4-GDD 20^5 b^1 50^1}.

{\boldmath $g = 24$}~
For $b \in \{33, 42, 51, 60\}$, use Theorem~\ref{thm:4-GDD (mg)^u} with a 4-GDD of type $8^5 (b/3)^1 20^1$.

For $b = 48$, use Theorem~\ref{thm:4-GDD (mg)^u} with a 4-GDD of type $6^5 12^1 15^1$ from \cite[Proposition IV.4.15]{Ge2007}.

For $b \in \{27, 30, 36, 39, 45, 54, 57\}$, see Lemma~\ref{lem:4-GDD 24^5 b^1 60^1}.

{\boldmath $g = 30$}~
For $b \in \{48, 66\}$, use Theorem~\ref{thm:4-GDD (mg)^u} with a 4-GDD of type $10^5 (b/3)^1 25^1$.

For $b  = 60$, use Theorem~\ref{thm:4-GDD (mg)^u} with a 4-GDD of type $6^5 12^1 15^1$ from \cite[Proposition IV.4.15]{Ge2007}.

For $b \in \{36, 42, 54, 72\}$, see Lemma~\ref{lem:4-GDD 30^5 b^1 75^1}.

{\boldmath $g = 40$}~
For $b \in \{55, 70, 85, 100\}$, use Theorem~\ref{thm:4-GDD (mg)^u} with a 4-GDD of type $8^5 (b/5)^1 20^1$.

For $b \in \{64, 88\}$, use Theorem~\ref{thm:4-GDD (mg)^u} with a 4-GDD of type $10^5 (b/4)^1 25^1$.

For $b \in \{43, 46, 49, 52, 58, 61, 67, 73, 76, 79, 82, 91, 94, 97\}$, see Lemma~\ref{lem:4-GDD 40^5 b^1 100^1}.

{\boldmath $g = 60$}~
For $b \in \{69, 78, 87, 96, 105, 114, 123, 132, 141, 150\}$, use Theorem~\ref{thm:4-GDD (mg)^u} with a 4-GDD of type $20^5 (b/3)^1 50^1$.

For $b \in \{75, 90, 120, 135\}$, use Theorem~\ref{thm:4-GDD (mg)^u} with a 4-GDD of type $12^5 (b/5)^1 30^1$.

For $b \in \{63, 66, 72, 81, 84, 93, 99, 102, 108, 111, 117, 126, 129, 138, 144, 147\}$, see Lemma~\ref{lem:4-GDD 60^5 b^1 150^1}.

{\boldmath $g = 80$}~
For even $b \le 176$, use Theorem~\ref{thm:(am)^d b^1 (adm/2)^1} with $\epsilon = 6$, $d = 5$, $m = 8$, $\alpha = 10$ and
4-GDDs of type $10^5 b^1 25^1$, $b \in \{10, 16, 22\}$.

For $b \in \{95, 125, 155, 185, 200\}$, use Theorem~\ref{thm:4-GDD (mg)^u} with a 4-GDD of type $16^5 (b/5)^1 40^1$.

For $b = 188$, use Theorem~\ref{thm:4-GDD (mg)^u} with a 4-GDD of type $20^5 47^1 50^1$.

For $b \in$ \{83, 89, 101, 107, 113, 119, 131, 137, 143, 149, 161, 167, 173, 179, 182, 191, 194, 197\},
see Lemma~\ref{lem:4-GDD 80^5 b^1 200^1}.

{\boldmath $g = 120$}~
Use Theorem~\ref{thm:(am)^d b^1 (adm/2)^1 hole} with
an ITD$(7, 30; 5)$ from \cite[Table III.4.14]{AbelColbournDinitz2007b} and
4-GDDs of types $4^5 \beta^1 10^1$ and $20^5 \gamma^1 50^1$.
~\adfQED
%
%%%%%%%%%%%%%%%%%%%%%%%%%%%%%%%%%%%%%%%%%%%%%%%%%%%%%%%%%%%%%%%%%%%%%%%%%%%%%%%%%%%%%%%%%%
%%%%%%%%%%%%%%%%%%%%%%%%%%%%%%%%%%%%%%%%%%%%%%%%%%%%%%%%%%%%%%%%%%%%%%%%%%%%%%%%%%%%%%%%%%
\begin{theorem}
\label{thm:4-GDD (4a)^5 b^1 (10a)^1 existence}
Suppose $g \equiv 0 \adfmod{4}$.
A 4-GDD of type $g^5 b^1 (5g/2)^1$ exists if and only if $g \ge 4$, $b \equiv g \adfmod{3}$ and $g \le b \le 5g/2$.
\end{theorem}
{\bf Proof}~ Clearly, we must have $g \ge 4$ and $b \equiv g \adfmod{3}$. Designate points in the 5 groups of size $g$ as type $A$,
points in the group of size $b$ as type $B$, and
points in the group of size $5 g/2$ as type $C$.
Then the 4-GDD has $5g (b + 3 g)/4$ blocks, broken down by type as follows:
\begin{center}
\begin{tabular}{lr@{~~~~~~}lr}
type $AAAA$ & $5 (b - g) g / 12$,     & type $AAAB$ & 0, \\
type $AAAC$ & $5 (5g - 2 b) g / 6$,   & type $AABC$ & $5 b g/2$,
\end{tabular}
\end{center}
and hence $g \le b \le 5g/2$. Therefore the stated conditions are necessary.

For sufficiency, there exist five MOLS of side $m$ for all $m \ge 7$ except possibly
$$m \in \{10, 14, 15, 20, 22, 26, 30, 34, 38, 46\};$$
see \cite{Abel2015} for $m \in \{18, 60\}$, \cite[Table III.3.88]{AbelColbournDinitz2007} for the others.
Therefore, by Theorem~\ref{thm:(am)^d b^1 (adm/2)^1} with $\epsilon = 3$ and $\alpha = 4$, together with Theorem~\ref{thm:(4a)^5 b^1 (10a)^1 small},
for all $m \ge 7$ except possibly
\begin{equation}
\label{eqn:4-GDD (4a)^5 b^1 (10a)^1 exceptions}
4m \in \{40, 56, 60, 80, 88, 104, 120, 136, 152, 184\},
\end{equation}
there exist 4-GDDs of type $(4m)^5 b^1 (10m)^1$ for $b \in \{4m, 4m+3, \dots, 10m\}$.
Moreover, by Theorem~\ref{thm:(4a)^5 b^1 (10a)^1 small} and Theorem~\ref{thm:(am)^d b^1 (adm/2)^1} with $\epsilon = 3$ and $\alpha = 8$,
we can eliminate 56, 88, 104, 136, 152 and 184 from (\ref{eqn:4-GDD (4a)^5 b^1 (10a)^1 exceptions}).
Theorem~\ref{thm:(4a)^5 b^1 (10a)^1 small} deals directly with 40, 60, 80 and 120.
\adfQED
%
%%%%%%%%%%%%%%%%%%%%%%%%%%%%%%%%%%%%%%%%%%%%%%%%%%%%%%%%%%%%%%%%%%%%%%%%%%%%%%%%%%%%%%%%%%
%%%%%%%%%%%%%%%%%%%%%%%%%%%%%%%%%%%%%%%%%%%%%%%%%%%%%%%%%%%%%%%%%%%%%%%%%%%%%%%%%%%%%%%%%%
\begin{theorem}
\label{thm:(6a+3)^6 b^1 (18a+9)^1 small}
There exists a 4-GDD of type $g^6 b^1 (3g)^1$ for
$g \in$ \{$3$, $9$, $15$, $45$\} and $b \in \{g, g + 6, \dots, 3g\}$.
\end{theorem}
{\bf Proof}~ By Theorem~\ref{thm:4-GDD g^u m^1 existence}, we may assume $b \ge g + 6$.
We deal with each value of $g$ in turn.

{\boldmath $g = 3$}~
For $b = 9$, see \cite[Proposition IV.4.15]{Ge2007}.

{\boldmath $g = 9$}~
For $b = 15$, see  \cite[Lemma 2.1]{Forbes2019}.

For $b = 27$, use Theorem~\ref{thm:4-GDD (mg)^u} with a 4-GDD of type $3^6 9^1 9^1$.

For $b = 21$, see Lemma~\ref{lem:4-GDD 9^6 b^1 27^1}.

{\boldmath $g = 15$}~
For $b = 45$, use Theorem~\ref{thm:4-GDD (mg)^u} with a 4-GDD of type $3^6 9^1 9^1$.

For $b \in \{21, 27, 33, 39\}$, see Lemma~\ref{lem:4-GDD 15^6 b^1 45^1}.

{\boldmath $g = 45$}~
For $b \in \{63, 81, 99, 117, 135\}$, use Theorem~\ref{thm:4-GDD (mg)^u} with a 4-GDD of type $15^6 (b/3)^1 45^1$.

For $b \in \{75, 105\}$, use Theorem~\ref{thm:4-GDD (mg)^u} with a 4-GDD of type $9^6 (b/5)^1 27^1$.

For $b \in \{51, 57, 69, 87, 93, 111, 123, 129\}$, see Lemma~\ref{lem:4-GDD 45^6 b^1 135^1}.
~\adfQED
%
%%%%%%%%%%%%%%%%%%%%%%%%%%%%%%%%%%%%%%%%%%%%%%%%%%%%%%%%%%%%%%%%%%%%%%%%%%%%%%%%%%%%%%%%%%
%%%%%%%%%%%%%%%%%%%%%%%%%%%%%%%%%%%%%%%%%%%%%%%%%%%%%%%%%%%%%%%%%%%%%%%%%%%%%%%%%%%%%%%%%%
\begin{theorem}
\label{thm:4-GDD (6a+3)^6 b^1 (18a+9)^1 existence}
Suppose $g \equiv 3 \adfmod{6}$.
A 4-GDD of type $g^6 b^1 (3g)^1$ exists if and only if
$g \ge 3$, $b \equiv 3 \adfmod{6}$ and $3 \le b \le 3g$.
\end{theorem}
{\bf Proof}~ Necessity follows from an analysis of the block types similar to that in Theorem~\ref{thm:4-GDD (4a)^5 b^1 (10a)^1 existence}.
For sufficiency, see Theorem~\ref{thm:(6a+3)^6 b^1 (18a+9)^1 small} when $g \in \{3, 9, 15, 45\}$.
There exist 6 MOLS of side $q$ for all odd $q \ge 7$ except possibly
$$q \in \{15, 21, 33, 39, 51\};$$
see \cite{JaniszczakStaszewski2019} for 35, \cite[Table III.3.88]{AbelColbournDinitz2007} for the rest.
If $g \ge 21$ and $g \neq 45$, we can find $\alpha \in \{3, 9\}$ such that there exist 6 MOLS of side $m = g/\alpha$.
Hence we can use Theorem~\ref{thm:(am)^d b^1 (adm/2)^1} with $\epsilon = d = 6$.
\adfQED
%
%%%%%%%%%%%%%%%%%%%%%%%%%%%%%%%%%%%%%%%%%%%%%%%%%%%%%%%%%%%%%%%%%%%%%%%%%%%%%%%%%%%%%%%%%%
%%%%%%%%%%%%%%%%%%%%%%%%%%%%%%%%%%%%%%%%%%%%%%%%%%%%%%%%%%%%%%%%%%%%%%%%%%%%%%%%%%%%%%%%%%
\begin{theorem}
\label{thm:(6a)^6 b^1 (18a)^1 small}
There exists a 4-GDD of type $g^6 b^1 (3g)^1$ for
$$g \in \{6, 12, 18, 24, 30, 36, 60, 72, 90, 120, 180, 360\}$$
and $b \in \{g, g + 3, \dots, 3g\}$.
\end{theorem}
{\bf Proof}~ By Theorem~\ref{thm:4-GDD g^u m^1 existence}, we may assume $b \ge g + 3$.
We deal with each value of $g$ in turn.
We refer to
Theorems~\ref{thm:4-GDD (3a)^4 (3b)^1 (6a)^1 existence} and
\ref{thm:4-GDD (4a)^5 b^1 (10a)^1 existence} as well as earlier results of this theorem
for the existence of the 4-GDDs used in Theorems~\ref{thm:(6am)^d b^1 (3adm)^1} and \ref{thm:(am)^d b^1 (adm/2)^1 hole}.

{\boldmath $g = 6$}~
For $b = 9$, see \cite[Proposition IV.4.15]{Ge2007}.

For $b \in \{12, 15, 18\}$, see Lemma~\ref{lem:4-GDD 6^6 b^1 18^1}.

{\boldmath $g = 12$}~
For $b \in \{15, 21\}$, see \cite[Proposition IV.4.15]{Ge2007}.

For $b \in \{18, 24, 27, 30, 33, 36\}$, see Lemma~\ref{lem:4-GDD 12^6 b^1 36^1}.

{\boldmath $g = 18$}~
For $b = 21$, see \cite[Lemma 2.1]{Forbes2019}.

For $b \in \{27, 36, 45, 54\}$, use Theorem~\ref{thm:4-GDD (mg)^u} with a 4-GDD of type $6^6 (b/3)^1 18^1$.

For $b \in \{24, 30, 33, 39, 42, 48, 51\}$, see Lemma~\ref{lem:4-GDD 18^6 b^1 54^1}.

{\boldmath $g = 24$}~
For $b = 27$, see \cite[Lemma 2.1]{Forbes2019}.

For $b \in \{36, 48, 60, 72\}$, use Theorem~\ref{thm:4-GDD (mg)^u} with a 4-GDD of type $6^6 (b/4)^1 18^1$.

For $b \in \{30, 33, 39, 42, 45, 51, 54, 57, 63, 66, 69\}$, see Lemma~\ref{lem:4-GDD 24^6 b^1 72^1}.

{\boldmath $g = 30$}~
For $b \le 60$, use Theorem~\ref{thm:(6am)^d b^1 (3adm)^1} with $m = 5$, $d = 6$, $\alpha = 1$ and
4-GDDs of types $6^4 \beta^1 12^1$ and $6^6 \beta^1 18^1$, $\beta \in \{6, 9, 12\}$.

For $b \in \{75, 90\}$, use Theorem~\ref{thm:4-GDD (mg)^u} with a 4-GDD of type $6^6 (b/5)^1 18^1$.

For $b \in \{63, 66, 69, 72, 78, 81, 84, 87\}$, see Lemma~\ref{lem:4-GDD 30^6 b^1 90^1}.

{\boldmath $g = 36$}~
For $b \in \{45, 54, 63, 72, 81, 90, 99, 108\}$, use Theorem~\ref{thm:4-GDD (mg)^u} with a 4-GDD of type $12^6 (b/3)^1 36^1$.

For $b \in \{39, 42, 48, 51, 57, 60, 66, 69, 75, 78, 84, 87, 93, 96, 102, 105\}$, see Lemma \ref{lem:4-GDD 36^6 b^1 108^1}.

{\boldmath $g = 60$}~
For $b \le 120$, use Theorem~\ref{thm:(6am)^d b^1 (3adm)^1} with $m = 5$, $d = 6$, $\alpha = 2$ and
4-GDDs of types $12^4 \beta^1 24^1$ and $12^6 \beta^1 36^1$, $\beta \in \{12, 15, \dots, 24\}$.

For $b \in \{135, 150, 165, 180\}$, use Theorem~\ref{thm:4-GDD (mg)^u} with a 4-GDD of type $12^6 (b/5)^1 36^1$.

For $b \in \{132, 156\}$, use Theorem~\ref{thm:4-GDD (mg)^u} with a 4-GDD of type $15^6 (b/4)^1 45^1$.

For $b \in \{123, 126, 129, 138, 141, 144, 147, 153, 159, 162, 168, 171, 174, 177\}$, see Lemma~\ref{lem:4-GDD 60^6 b^1 180^1}.

{\boldmath $g = 72$}~
For $b \le 180$, use Theorem~\ref{thm:(6am)^d b^1 (3adm)^1} with $m = 6$, $d = 6$, $\alpha = 2$ and
4-GDDs of types $12^5 \beta^1 30^1$ and $12^6 \beta^1 36^1$, $\beta \in \{12, 15, \dots, 30\}$.

For even $b$, use Theorem~\ref{thm:(am)^d b^1 (adm/2)^1} with $\epsilon = d = 6$, $m = 24$ and
4-GDDs of types $3^6 3^1 9^1$ and $3^6 9^1 9^1$ from Theorem~\ref{thm:4-GDD (6a+3)^6 b^1 (18a+9)^1 existence}.

For $b \in \{189, 207\}$, use Theorem~\ref{thm:4-GDD (mg)^u} with a 4-GDD of type $24^6 (b/3)^1 72^1$.

For $b \in \{183, 195, 201, 213\}$, see Lemma~\ref{lem:4-GDD 72^6 b^1 216^1}.

{\boldmath $g = 90$}~
For $b \le 180$, use Theorem~\ref{thm:(6am)^d b^1 (3adm)^1} with $m = 5$, $d = 6$, $\alpha = 3$ and
4-GDDs of types $18^4 \beta^1 36^1$ and $18^6 \beta^1 54^1$, $\beta \in \{18, 21, \dots, 36\}$.

For $b \in \{189, 198, 207, 216, 225, 234, 243, 252, 261, 270\}$, use Theorem~\ref{thm:4-GDD (mg)^u} with a 4-GDD of type $30^6 (b/3)^1 90^1$.

For $b \in \{195, 210, 240, 255\}$, use Theorem~\ref{thm:4-GDD (mg)^u} with a 4-GDD of type $18^6 (b/5)^1 54^1$.

For $b \in \{183, 186, 192, 201, 204, 213, 219, 222, 228, 231, 237, 246, 249, 258, 264, 267\}$, see Lemma~\ref{lem:4-GDD 90^6 b^1 270^1}.

{\boldmath $g = 120$}~
For $b \le 240$, use Theorem~\ref{thm:(6am)^d b^1 (3adm)^1} with $m = 5$, $d = 6$, $\alpha = 4$ and
4-GDDs of types $24^4 \beta^1 48^1$ and $24^6 \beta^1 72^1$, $\beta \in \{24, 27, \dots, 48\}$.

For even $b$, use Theorem~\ref{thm:(am)^d b^1 (adm/2)^1} with $m = 40$, $d = 6$, $\alpha = 3$ and
4-GDDs of types $3^6 \beta^1 9^1$, $\beta \in \{3, 9\}$.
By \cite[Table III.3.88]{AbelColbournDinitz2007}, there exist 6 MOLS of side $40$.

For $b \in \{255, 285, 315, 345\}$, use Theorem~\ref{thm:4-GDD (mg)^u} with a 4-GDD of type $24^6 (b/5)^1 72^1$.

For $b \in$ \{243, 249, 261, 267, 273, 279, 291, 297, 303, 309, 321, 327, 333, 339, 351, 357\}, see Lemma~\ref{lem:4-GDD 120^6 b^1 360^1}.

{\boldmath $g = 180$}~
For $b \le 360$, use Theorem~\ref{thm:(6am)^d b^1 (3adm)^1} with $m = 5$, $d = 6$, $\alpha = 6$ and
4-GDDs of types $36^4 \beta^1 72^1$ and $36^6 \beta^1 108^1$, $\beta \in \{36, 39, \dots, 72\}$.

For even $b$, use Theorem~\ref{thm:(am)^d b^1 (adm/2)^1 hole} with $\epsilon = 6$,
an ITD$(8, 60; 4)$ from \cite[Table III.4.14]{AbelColbournDinitz2007b} and
4-GDDs of types $3^6 \beta^1 9^1$, $\beta \in \{3,9\}$, and $12^6 \gamma^1 36^1$, $\gamma \in \{12, 15, \dots, 36\}$.

For $b \in \{369, 387, 405, 423, 441, 459, 477, 495, 513, 531\}$, use Theorem~\ref{thm:4-GDD (mg)^u} with a 4-GDD of type $60^6 (b/3)^1 180^1$.

For $b \in \{375, 435, 465, 525\}$, use Theorem~\ref{thm:4-GDD (mg)^u} with a 4-GDD of type $36^6 (b/5)^1 108^1$.

For $b \in$ \{363, 381, 393, 399, 411, 417, 429, 447, 453, 471, 483, 489, 501, 507, 519, 537\}, see Lemma~\ref{lem:4-GDD 180^6 b^1 540^1}.

{\boldmath $g = 360$}~
Use Theorem~\ref{thm:(am)^d b^1 (adm/2)^1 hole} with
an ITD$(8, 60; 4)$ from \cite[Table III.4.14]{AbelColbournDinitz2007b} and
4-GDDs of types $6^6 \beta^1 18^1$, $\beta \in \{6, 9, 12, 15, 18\}$, and $24^6 \gamma^1 72^1$, $\gamma \in \{24, 27, \dots, 72\}$.
~\adfQED
%
%%%%%%%%%%%%%%%%%%%%%%%%%%%%%%%%%%%%%%%%%%%%%%%%%%%%%%%%%%%%%%%%%%%%%%%%%%%%%%%%%%%%%%%%%%
%%%%%%%%%%%%%%%%%%%%%%%%%%%%%%%%%%%%%%%%%%%%%%%%%%%%%%%%%%%%%%%%%%%%%%%%%%%%%%%%%%%%%%%%%%
\begin{theorem}
\label{thm:4-GDD (6a)^6 b^1 (18a)^1 existence}
Suppose $g \equiv 0 \adfmod{6}$.
A 4-GDD of type $g^6 b^1 (3g)^1$ exists if and only if $g \ge 6$, $b \equiv 0 \adfmod{3}$ and $g \le b \le 3g$.
\end{theorem}
{\bf Proof}~ Necessity follows from an analysis of the block types similar to that in Theorem~\ref{thm:4-GDD (4a)^5 b^1 (10a)^1 existence}.
For sufficiency, see Theorem~\ref{thm:(6a)^6 b^1 (18a)^1 small} if $g \le 36$.
By \cite[Table III.3.88]{AbelColbournDinitz2007}, and \cite{JaniszczakStaszewski2019} for $q = 35$, there exist 6 MOLS of side $q$ for all $q \ge 7$ except possibly
\begin{align*}
q \in\; & \{10, 12, 14, 15, 18, 20, 21, 22, 26, 28, 30, 33, 34, 38, 39, \\
          &~~ 42, 44, 46, 51, 52, 54, 58, 60, 62, 66, 68, 74\}.
\end{align*}
Let $X_6 = \{60, 72, 90, 120, 180, 360\}$.
If $g \ge 42$ and $g \not\in X_6$, we can find $\delta \le 6$ such that there exist 6 MOLS of side $g/(6\delta)$.
Hence we can use Theorem~\ref{thm:(am)^d b^1 (adm/2)^1} with $\epsilon = 3$, $d = 6$, $m = g/(6\delta)$, $\alpha = 6 \delta$
and 4-GDDs of type $(6\delta)^6 \beta^1 (18\delta)^1$ from Theorem~\ref{thm:(6a)^6 b^1 (18a)^1 small}.
For $g \in X_6$, see Theorem~\ref{thm:(6a)^6 b^1 (18a)^1 small}.
\adfQED
%
%%%%%%%%%%%%%%%%%%%%%%%%%%%%%%%%%%%%%%%%%%%%%%%%%%%%%%%%%%%%%%%%%%%%%%%%%%%%%%%%%%%%%%%%%%
%%%%%%%%%%%%%%%%%%%%%%%%%%%%%%%%%%%%%%%%%%%%%%%%%%%%%%%%%%%%%%%%%%%%%%%%%%%%%%%%%%%%%%%%%%
\begin{theorem}
\label{thm:(12a)^7 b^1 (42a)^1 small}
There exists a 4-GDD of type $g^7 b^1 (7g/2)^1$ for
$g \in$ \{$12$, $24$, $36$, $48$, $60$, $72$, $84$\} and $b \in \{g, g + 3, \dots, 7g/2\}$.
\end{theorem}
{\bf Proof}~ By Theorem~\ref{thm:4-GDD g^u m^1 existence}, we may assume $b \ge g + 3$.
We deal with each value of $g$ in turn.

{\boldmath $g = 12$}~
For $b = 15$, see \cite[Lemma 2.1]{Forbes2019}.

For $b \in \{18, 21, \dots, 42\}$, see Lemma~\ref{lem:4-GDD 12^7 b^1 42^1}.

{\boldmath $g = 24$}~
For $b = 27$, see \cite[Lemma 2.1]{Forbes2019}.

For $b \in \{48, 72\}$,
use Theorem~\ref{thm:4-GDD (mg)^u} with a 4-GDD of type $6^7 (b/4)^1 21^1$ from Lemma~\ref{lem:4-GDD 6^7 b^1 21^1}.

For $b \in \{30, 33, \dots, 84\} \setminus \{48, 72\}$, see Lemma~\ref{lem:4-GDD 24^7 b^1 84^1}.

{\boldmath $g = 36$}~
For $b = 39$, see \cite[Lemma 2.1]{Forbes2019}.

For $b \in \{45, 54, 63, 72, 81, 90, 99, 108, 117, 126\}$,
use Theorem~\ref{thm:4-GDD (mg)^u} with a 4-GDD of type $12^7 (b/3)^1 42^1$ (from this theorem).

For $b \in$ \{42, 48, 51, 57, 60, 66, 69, 75, 78, 84, 87, 93, 96, 102, 105, 111, 114, 120, 123\}, see Lemma~\ref{lem:4-GDD 36^7 b^1 126^1}.

{\boldmath $g = 48$}~
For $b = 51$, see \cite[Lemma 2.1]{Forbes2019}.

For $b \in \{60, 72, 84, 96, 108, 120, 132, 144, 156, 168\}$, use Theorem~\ref{thm:4-GDD (mg)^u} with a 4-GDD of type $12^7 (b/4)^1 42^1$.

For $b \in$ \{54, 57, 63, 66, 69, 75, 78, 81, 87, 90, 93, 99, 102, 105, 111, 114, 117, 123, 126, 129, 135, 138, 141, 147, 150, 153, 159, 162, 165\}, see Lemma~\ref{lem:4-GDD 48^7 b^1 168^1}.

{\boldmath $g = 60$}~
For $b \le 120$, use Theorem~\ref{thm:(6am)^d b^1 (3adm)^1} with $d = 7$, $m = 5$, $\alpha = 2$ and
4-GDDs of types $12^4 \beta^1 24^1$ from Theorem~\ref{thm:4-GDD (3a)^4 (3b)^1 (6a)^1 existence} and
$12^7 \beta^1 42^1$, $\beta \in \{12, 15, \dots, 24\}$.

For $b \in \{135, 150, 165, 180, 195, 210\}$, use Theorem~\ref{thm:4-GDD (mg)^u} with a 4-GDD of type $12^7 (b/5)^1 42^1$.

For $b \in$ \{123, 126, 129, 132, 138, 141, 144, 147, 153, 156, 159, 162, 168, 171, 174, 177, 183, 186, 189, 192, 198, 201, 204, 207\}, see Lemma~\ref{lem:4-GDD 60^7 b^1 210^1}.

{\boldmath $g = 72$}~
For $b \le 180$, use Theorem~\ref{thm:(6am)^d b^1 (3adm)^1} with $d = 7$, $m = 6$, $\alpha = 2$ and
4-GDDs of types $12^5 \beta^1 30^1$ from Theorem~\ref{thm:4-GDD (6a)^6 b^1 (18a)^1 existence} and
$12^7 \beta^1 42^1$, $\beta \in \{12, 15, \dots, 30\}$.

For $b \in \{189, 198, 207, 216, 225, 234, 243, 252\}$, use Theorem~\ref{thm:4-GDD (mg)^u} with a 4-GDD of type $24^7 (b/3)^1 84^1$.

For $b \in \{192, 240\}$, use Theorem~\ref{thm:4-GDD (mg)^u} with a 4-GDD of type $18^7 (b/4)^1 63^1$ from Lemma~\ref{lem:4-GDD 18^7 b^1 63^1}.

For $b \in$ \{183, 186, 195, 201, 204, 210, 213, 219, 222, 228, 231, 237, 246, 249\}, see Lemma~\ref{lem:4-GDD 72^7 b^1 252^1}.

{\boldmath $g = 84$}~
For $b \le 252$, use Theorem~\ref{thm:(6am)^d b^1 (3adm)^1} with $d = 7$, $m = 7$, $\alpha = 2$ and
4-GDDs of types $12^6 \beta^1 36^1$ from Theorem~\ref{thm:4-GDD (6a)^6 b^1 (18a)^1 existence} and
$12^7 \beta^1 42^1$, $\beta \in \{12, 15, \dots, 36\}$.

For $b \in \{273, 294\}$, use Theorem~\ref{thm:4-GDD (mg)^u} with a 4-GDD of type $12^7 (b/7)^1 42^1$.

For $b \in$ \{255, 258, 261, 264, 267, 270, 276, 279, 282, 285, 288, 291\}, see Lemma~\ref{lem:4-GDD 84^7 b^1 294^1}.
~\adfQED
%
%%%%%%%%%%%%%%%%%%%%%%%%%%%%%%%%%%%%%%%%%%%%%%%%%%%%%%%%%%%%%%%%%%%%%%%%%%%%%%%%%%%%%%%%%%
%%%%%%%%%%%%%%%%%%%%%%%%%%%%%%%%%%%%%%%%%%%%%%%%%%%%%%%%%%%%%%%%%%%%%%%%%%%%%%%%%%%%%%%%%%
\begin{theorem}
\label{thm:4-GDD (12a)^7 b^1 (42a)^1 existence}
Suppose $g \equiv 0 \adfmod{12}$.
A 4-GDD of type $g^7 b^1 (7g/2)^1$ exists if and only if $g \ge 12$, $b \equiv 0 \adfmod{3}$ and $g \le b \le 7g/2$,
except possibly for the following:
\begin{align*}
g \in &\{120, 180, 240, 360, 420, 720, 840\} \textrm{~and~} 2g < b < 7g/2,\\
g \in &\{144, 1008\} \textrm{~and~} 5g/2 < b < 7g/2,\\
g \in &\{168, 252, 336, 504, 1512\} \textrm{~and~} 3g < b < 7g/2.
\end{align*}
\end{theorem}
{\bf Proof}~ Necessity follows from an analysis of the block types similar to that in Theorem~\ref{thm:4-GDD (4a)^5 b^1 (10a)^1 existence}.
For sufficiency, see Theorem~\ref{thm:(12a)^7 b^1 (42a)^1 small} if $g \le 84$.
By \cite[Table III.3.88]{AbelColbournDinitz2007}, and \cite{JaniszczakStaszewski2019} for $q = 63$,
there exist 7 MOLS of side $q$ for all $q \ge 8$ except possibly
$q \in$ \{10, 12, 14, 15, 18, 20, 21, 22, 26, 28, 30, 33,
34, 35, 38, 39, 42, 44, 45, 46, 50, 51, 52, 54, 55, 58, 60, 62, % 63,
66, 68, 69, 70, 74, 76, 77, 78, 84, 85, 86, 87, 90, 92, 93, 94, 95,
98, 102, 106, 108, 110, 111, 114, 116, 118, 119, 122, 123, 124, 126,
130, 132, 134, 138, 140, 142, 146, 148, 150, 156, 159, 162, 164, 166,
170, 172, 174, 175, 178, 180, 182, 183, 186, 188, 190, 194, 196, 198,
202, 204, 206, 212, 214, 218, 220, 222, 226, 228, 230, 234, 236, 238,
242, 244, 246, 250, 252, 258, 260, 274, 278, 282, 284, 286, 290, 291,
292, 294, 295, 306, 322, 326, 330, 335, 338, 340, 346, 348, 354, 358,
362, 366, 426, 430, 434, 436, 478, 482, 486, 490, 492, 494, 498, 506,
510, 566, 570\}.
Let
\begin{align*}
G_7 = {} & \{120, 144, 168, 180, 240, 252, 336, 360, 420, 504, 720, 840, \\
         &~~ 1008, 1080, 1320, 1512, 1584, 1680, 1872, 2040, 2280, 2448, \\
         &~~ 2640, 2856, 4080, 4176\},\\
X_7 = {} & \{120, 144, 168, 180, 240, 252, 336, 360, 420, 504, 720, 840, 1008, 1512\}.
% G_7/12 = \{10, 12, 14, 15, 20, 21, 28, 30, 35, 42, 60, 70, 84, 90, 110, 126, 132, 140, 156, 170, 190, 204, 220, 238, 340, 348\}
\end{align*}
If $g \ge 96$ and $g \not\in G_7$, we can find $\delta \le 7$ such that there exist 7 MOLS of side $g/(12\delta)$, and hence
we can use Theorem~\ref{thm:(am)^d b^1 (adm/2)^1} with $\epsilon = 3$, $d = 7$, $m = g/(12\delta)$ and $\alpha = 12 \delta$.

It remains to deal with $g \in G_7$.
We refer to Theorems~\ref{thm:4-GDD (3a)^4 (3b)^1 (6a)^1 existence},
\ref{thm:4-GDD (4a)^5 b^1 (10a)^1 existence},
\ref{thm:4-GDD (6a)^6 b^1 (18a)^1 existence},
\ref{thm:(12a)^7 b^1 (42a)^1 small} and earlier results of this theorem
for the existence of the 4-GDDs we need whenever we invoke
Theorems~\ref{thm:(6am)^d b^1 (3adm)^1} and \ref{thm:(am)^d b^1 (adm/2)^1 hole}.
We begin with $g \in G_7 \setminus X_7$,
referring to \cite[Table III.4.14]{AbelColbournDinitz2007b} for the existence of the required ITDs.
% Charlotte:GDD4-g^u-b^1-c^1-TeX-gen-A:HITS-fun:4.12
\adfDgap
%ADFvfyBlocksStart {6,6,6,6,6,6,18,18}
\noindent{\boldmath $ 6^{6} 18^{2} $}~
With the point set $Z_{72}$ partitioned into
 residue classes modulo $6$ for $\{0, 1, \dots, 35\}$, and
 residue classes modulo $2$ for $\{36, 37, \dots, 71\}$,
 the design is generated from
%ADFvfySectionEnd
\end{document}